\documentclass[a4paper]{article}
\usepackage{amsmath,amssymb,fancybox,longtable,graphics,amscd,multirow,theorem}
\usepackage[dvipdfm]{graphicx} \usepackage{caption,setspace}
\newtheorem{thm}{\textsc{Theorem}}[section]

\def\proof{\textsc{Proof. }}
\def\QED{$\Box$} 
\theorembodyfont{\rmfamily}

\theorembodyfont{\rmfamily}

\begin{document} 
\title{Polytope duality for families of $K3$ surfaces associated to singularities $Q_{16}$ and $S_{16}$}
\author{Makiko Mase}
\date{\empty}
\maketitle
\abstract{There are strange dual pairs of bimodal singularities that are not assigned an invertible projectivisation in \cite{EbelingPloog}. 
We study families of $K3$ surfaces associated to such pairs. } \\

%$K3$ and Enriques surfaces\quad Lattice polytopes (including relations with commutative algebra and algebraic geometry)\quad Singularity

\section{Introduction}
Ebeling and Takahashi \cite{EbelingTakahashi11} introduced a notion of strange duality for invertible polynomials by Berglund-H$\ddot{\textnormal{u}}$bsch mirror construction~\cite{BerglundHubsch}. 
Mase and Ueda \cite{MU} studied an extension of strange duality to polytope duality when bimodal singularities admit an invertible projectivisation, which are given in the study of Ebeling and Ploog \cite{EbelingPloog} of distinguished basis. 
In order to complete our list, we concern two pairs of strange duality : $Q_{16}$ and $S_{16}$ that are defined by 
\[
Q_{16} :  f_{Q_{16}} = x^4z+y^3+xz^2, \qquad
S_{16} :  f_{S_{16}} = x^4y+xz^2+y^2z.   
\]
The matrices of exponents of the defining polynomials are respectively given by 
\[
A_{f_{Q_{16}}} = \begin{pmatrix} 4 & 0 & 1 \\ 0 & 3 & 0 \\ 1 & 0 & 2\end{pmatrix}, \qquad
A_{f_{S_{16}}} = \begin{pmatrix} 4 & 1 & 0 \\ 1 & 0 & 2 \\ 0 & 2 & 1\end{pmatrix},
\]
which are both symmetric. 
Thus, the singularities are self-dual. 

In \cite{EbelingPloog} a projectivisations of these singularities are determined as follows: 
\[
Q_{16} :  F_{Q_{16}}=X^4Z+Y^3+XZ^2+W^6Z+W^7Y, 
\]
\[
S_{16} : F_{S_{16}}= X^4Y+XZ^2+Y^2Z+W^5Z+W^6Y. 
\]
Note that these polynomials are not invertible. 
Recall that a polynomial $F$ is {\it invertible} if the matrix of exponents of $F$ is an invertible matrix. 

The polynomial $F_{Q_{16}}$ is an anticanonical section of the weighted projective space $\mathbb{P}(2,3,7,9)$, and $F_{S_{16}}$ of $\mathbb{P}(2,3,5,7)$. 
It is known that both spaces $\mathbb{P}(2,3,7,9)$, and $\mathbb{P}(2,3,5,7)$ are Fano $3$-folds as is classified by Yonemura \cite{Yonemura}, thus general anticanonical sections are $K3$ surfaces with at most Gorenstein singularities as is shown by Batyrev \cite{BatyrevMirror}. 
It was concluded in \cite{MU} that all strange-dual pairs for bimodal singularities admitting an invertible projectivisation extend to polytope duality. 
Even though the singularities $Q_{16}$ and $S_{16}$ are not assigned such a projectivisation by \cite{EbelingPloog} thus have been omitted so far, it still makes sense to consider polytope duality of families of $K3$ surfaces. 

%Motivated by the homological mirror symmetry conjecture, we consider a pair of a polynomial $F$ and the {\it group of maximal diagonal symmetries} $G_{\rm max}(F)$ of $F$ defined as follows. 
%For a polynomial $F$ in $\mathbb{P}(w_0, w_1, w_2, w_3)$, define $G_{\rm max}(F)$ by 
%\[G_{\rm max}(F) = \{ (\alpha,\, \beta,\, \gamma,\, \delta)\in(\mathbb{C}^*)^4\, |\, F(\alpha W,\, \beta X,\, \gamma Y,\, \delta Z) = F(W,\, X,\, Y,\, Z)\}, \]
%where $(W,\, X,\, Y,\, Z)$ is a coordinate system of $\mathbb{P}(w_0, w_1, w_2, w_3)$. 
%For $G:=G_{\rm max}(F)$, denote by $\Delta_{a,G}$ and by $\Delta_{F, G}$ the convex hull of vertices that are corresponding to the anticanonical monomials in $\mathbb{P}(w_0, w_1, w_2, w_3)$ and monomials appearing in $F$ that are invariant by the action by $G$, respectively. 
%We call the polytope $\Delta_{F, G}$ the {\it Newton polytope of $F$}. 

In this article, we pose the following problem. 

\noindent
{\sc Problem. } 
For the polytope $F$ which is $F_{Q_{16}}$ or $F_{S_{16}}$, let $\Delta$ be a reflexive polytope such that $\Delta_F\subset\Delta\subset\Delta_w$, where $w$ is a weight system $(2,3,7,9;21)$, or $(2,3,5,7;17)$, respectively. 
Set $\Delta':=\Delta^*$ the polar dual polytope of $\Delta$. 
Determine whether or not the strange duality for the singularity extends to the polytope duality between subfamilies of $K3$ surfaces in $\mathbb{P}(w)$; 
equivalently, determine whether or not the polytope $\Delta'$ is a subpolytope of the polytope $\Delta_w$. 

The main theorem of the article is a negative answer to this question stated below. 

\noindent
{\bf Theorem \ref{MainThm}}\quad 
{\it 
Let $w$ be a weight system $(2,3,7,9)$ (resp. $(2,3,5,7)$) and $\Delta$ be a reflexive polytope such that $\Delta_F\subset\Delta\subset\Delta_w$, where $F$ is $F_{Q_{16}}$ (resp. $F_{S_{16}}$). 
The polar dual polytope $\Delta'$ of $\Delta$ is not a subpolytope of the polytope $\Delta_w$. 
In particular, the strange duality for the singularity $Q_{16}$ (resp, $S_{16}$) does not extend to the polytope duality between subfamilies of $K3$ surfaces in $\mathbb{P}(w)$. 
} \\

We recall Yonemura's list of weight systems and toric geometry, and then explain the polytope duality in section 2. 
The main theorem will be proved in section 3. 
We give a conclusion on a relation between the strange duality for bimodal singularities and the polytope duality in section 4. \\

\begin{ackn}
\textnormal{
The author thanks to Professor Wolfgang Ebeling for his suggestion of the study and discussions. }
\end{ackn}

\section{Preliminary}
A {\it $K3$ surface} is a compact complex non-singular $2$-dimensional algebraic variety with trivial canonical bundle and irregularity zero. 
It is known that a $K3$ surface with at most Gorenstein singularities, that is, rational double points, is birationally equivalent to a $K3$ surface due to the existence of a crepant resolution. 

A quadruple of positive integers $w=(w_0,\, w_1,\, w_2,\, w_3)$ is {\it well-posed} if $w_0\leq  w_1\leq w_2\leq w_3$, and any three $w_i,\, w_j,\, w_k$ out of four are prime. 
Set $d:=w_0+w_1+w_2+w_3$. 
Let $\mathbb{C}[W,\, X,\, Y,\, Z]$ be the ring of polynomials graded by the weights $\deg W=w_0,\, \deg X=w_1,\,\deg Y=w_2$, and $\deg Z=w_3$, with which we define the {\it weighted projective space} with {\it weight system $(w_0,\, w_1,\, w_2,\, w_3;\, d)$} by 
\[
\mathbb{P}(w)=\mathbb{P}(w_0, w_1, w_2, w_3):={\rm Proj}\,\mathbb{C}[W,\, X,\, Y,\, Z]. 
\]
It is known that weighted projective spaces are toric $3$-folds: we denote by $\Delta_w=\Delta_{(w_0, w_1, w_2, w_3)}$ the $3$-dimensional polytope associated to $\mathbb{P}(w)$. 

The anticanonical divisor of $\mathbb{P}(w)$ is isomorphic to $\mathcal{O}(-d)$ (see \cite{DolgachevWP}). 
Thus the global sections of it are polynomials of weighted degree $d$, which we simply call {\it anticanonical sections}. 
Indeed, let $M_w$ be a lattice of rank $3$ defined by 
\[
M_w:=\{ (i,j,k,l)\in\mathbb{Z}^4\, |\, w_0i+w_1j+w_2k+w_3l=0\}
\]
with a basis, $\Delta_{w}$ is embedded into $\mathbb{R}^3$, and a monomial $W^{i+1}X^{j+1}Y^{k+1}Z^{l+1}$ of weighted degree $d$ is corresponding to an element $(i, j, k, l)$ of $M_w$. 

The weighted projective spaces that are Fano, namely, the anticanonical divisor is ample, are classified by Yonemura into $95$ classes \cite{Yonemura}. 

Let $M$ be a lattice of rank $3$, and $N:={\rm Hom}_{\mathbb{Z}}(M,\,\mathbb{Z})$ be its dual with a natural pairing $\langle\, , \,\rangle:N\times M\to\mathbb{Z}$ and $\langle\, ,\,\rangle_{\mathbb{R}}$ is the extension to $\mathbb{R}$-coefficients, and  $\Delta$ be a $3$-dimensional convex hull of finite number of points in $M\otimes_{\mathbb{Z}}\mathbb{R}$, which we simply call a {\it polytope}. 
Define the {\it polar dual} polytope $\Delta^*$ of $\Delta$ by 
\[
\Delta^* := \{ y\in N\otimes_{\mathbb{Z}}\mathbb{R} \, |\, \langle y , x\rangle_{\mathbb{R}}\geq -1\: \textnormal{for all}\: x\in\Delta\}. 
\]
A polytope $\Delta$ with all vertices being integral points is {\it reflexive} if $\Delta$ contains the only integral points in its interior and the polar dual $\Delta^*$ has also all vertices integral. 
In general, if a polytope $\Delta$ is reflexive, the associated projective space $\mathbb{P}_{\!\Delta}$ is a Fano $3$-fold, and its general anticanonical sections are $K3$ surfaces with at most Gorenstein singularities (see \cite{BatyrevMirror}). 

For non-degenerate isolated singularities $(f=0)$ and $(f'=0)$ in $\mathbb{C}^3$ with projectivisations $F$ in $\mathbb{P}(w)=\mathbb{P}(w_0, w_1, w_2, w_3)$ and $F'$ in $\mathbb{P}(w')=\mathbb{P}(w'_0, w'_1, w'_2, w'_3)$ as anticanonical sections of Fano weighted projective spaces, families $\mathcal{F}_\Delta$ and $\mathcal{F}_{\Delta'}$ of $K3$ surfaces with at most Gorenstein singularities associated to reflexive polytopes $\Delta$ and $\Delta'$ ({\it c.f.} \cite{BatyrevMirror}) are said {\it polytope dual} if the following relations hold: 
\[
\Delta_F\subset\Delta\subset\Delta_{(w)},\quad
\Delta_{F'}\subset\Delta'\subset\Delta_{(w')},\quad \textnormal{and}\quad
\Delta^*\simeq\Delta'. 
\]

\section{Main Result}
%Recall that for a polynomial $F$ in $\mathbb{P}(w_0, w_1, w_2, w_3)$, the group of maximal diagonal symmetries $G_{\rm max}(F)$ is defined by 
%\[G_{\rm max}(F) = \{ (\alpha,\, \beta,\, \gamma,\, \delta)\in(\mathbb{C}^*)^4\, |\, F(\alpha W,\, \beta X,\, \gamma Y,\, \delta Z) = F(W,\, X,\, Y,\, Z)\}. \]

Recall our main theorem that is proved in this section. 

\begin{thm}\label{MainThm}
Let $w$ be a weight system $(2,3,7,9)$ (resp. $(2,3,5,7)$) and $\Delta$ be a reflexive polytope such that $\Delta_F\subset\Delta\subset\Delta_w$, where $F$ is $F_{Q_{16}}$ (resp. $F_{S_{16}}$). 
The polar dual polytope $\Delta'$ of $\Delta$ is not a subpolytope of the polytope $\Delta_w$. 
In particular, the strange duality for the singularity $Q_{16}$ (resp, $S_{16}$) does not extend to the polytope duality between subfamilies of $K3$ surfaces in $\mathbb{P}(w)$. 
\end{thm}
\proof

\noindent
{\bf Singularity of type }$\mathbf{Q_{16}}. $
The singularity is defined by a polynomial $f = x^4z+y^3+xz^2$, and take a projectivisation $F=X^4Z+Y^3+XZ^2+W^6Z+W^7Y$ in accordance of \cite{EbelingPloog}. 
%By definition, the group of maximal diagonal symmetries is given by 
%\begin{eqnarray*}G:=G_{\rm max}(F) & = & \{ (\alpha,\, \beta,\, \gamma,\, \delta)\in(\mathbb{C}^*)^4\, |\, \beta^4\delta = \gamma^3=\beta\delta^2=\alpha^6\delta=\alpha^7\gamma=1\} \\  & =& \left\langle \frac{1}{21}(2,3,7,9),\, \frac{1}{3}(2,0,1,0)\right\rangle. \end{eqnarray*}
%Here, $\frac{1}{n}(a_0,\, a_1,\, a_2,\, a_3)$ denotes the action 
%\[W^{x_0}X^{x_1}Y^{x_2}Z^{x_3} \mapsto (e^{\frac{2\pi i}{n}a_0}W)^{x_0}(e^{\frac{2\pi i}{n}a_1}X)^{x_1}(e^{\frac{2\pi i}{n}a_2}Y)^{x_2}(e^{\frac{2\pi i}{n}a_3}Z)^{x_3}. \]
Let $M$ be a lattice of rank $3$ defined by 
\[
M:=\left\{ (i,j,k,l)\in\mathbb{Z}^4 \, | \, 2i+ 3j+ 7k+ 9l = 0\right\}. 
\]

%It is easily seen that all monomials of weighted degree $21$ in $\mathbb{P}(2,3,7,9)$ are invariant by the action of $G$. 
%Thus, in this case, $\Delta_{(2,3,7,9), G}$ coincides with $\Delta_{(2,3,7,9)}$. 

By taking a basis $\{ e_1,\, e_2,\, e_3\}$ of $M$ by 
\[
e_1 = (8, 0, -1, -1),\, 
e_2 = (6, -1, 0, -1),\, 
e_3 = (5, -1, -1, 0), 
\]
the polytope $\Delta_{(2,3,7,9)}$ is given by a convex hull of vertices 
\[
(1, 0, 0),\, (0, 1, 0),\, (0, 0, 1),\, (0, -1, 1),\, (-1, 2, -1),\, (4, -3, -3), 
\]
respectively corresponding to monomials 
\[
W^9X,\, W^7Y,\, W^6Z,\, XZ^2,\, Y^3,\, X^7. 
\]

The Newton polytope $\Delta_F$ is a convex hull of vertices 
\[
(0, 1, 0),\, (0, 0, 1),\, (0, -1, 1),\, (-1, 2, -1),\, (2, -2, -1), 
\]
where $(2, -2, -1)$ is corresponding to the monomial $X^4Z$. 
Since the polar dual of the face 
\[
{\rm Conv} \left\{ (0, 1, 0),\, (-1, 2, -1),\, (2, -2, -1) \right\}
\]
is a rational vertex $(-4/3, -1, 1/3)$, the Newton polytope is not reflexive. 

Any reflexive polytope $\Delta$ satisfying $\Delta_F\subset\Delta\subset\Delta_{(2,3,7,9)}$ is thus a convex hull of vertices
\[
(1, 0, 0),\, (0, 1, 0),\, (0, 0, 1),\, (0, -1, 1),\, (-1, 2, -1),\, (2, -2, -1),\, (n+1, -n, -n)
\]
with $n=1$. 

As long as there is an edge 
\[
\Gamma = {\rm Conv}\{ (0, -1, 1),\, (-1, 2, -1)\} 
\]
in $\Delta$, since the polar dual of $\Gamma$ is 
\[
\Gamma^* = {\rm Conv}\{ (8, 6, 5),\, (2, 0, -1)\}, 
\]
and thus there are $5$ lattice points in $\Gamma^*$, the polar dual polytope $\Delta':=\Delta^*$ of $\Delta$ should contain an edge with $5$ lattice points. 
However, by a direct check, there does not exist such an edge in or inside of the polytope $\Delta_{(2,3,5,7)}$. 
Thus, the polytope $\Delta'$ is not a subpolytope of $\Delta_{(2,3,5,7)}$. 
Therefore the assertion is verified. 

\noindent
{\bf Singularity of type }$\mathbf{S_{16}}. $
The singularity is defined by a polynomial $f = x^4y+xz^2+y^2z$, and take a projectivisation $F= X^4Y+XZ^2+Y^2Z+W^5Z+W^6Y$ in accordance of \cite{EbelingPloog}. 
%By definition, the group of maximal diagonal symmetries is given by 
%\begin{eqnarray*}G:=G_{\rm max}(F) & = & \{ (\alpha,\, \beta,\, \gamma,\, \delta)\in(\mathbb{C}^*)^4\, |\, \beta^4\gamma+\beta\delta^2+\gamma^2\delta+\alpha^5\delta+\alpha^6\gamma=1\} \\ & = & \left\langle \frac{1}{17}(2,3,5,7)\right\rangle. \end{eqnarray*}
%Here, $\frac{1}{n}(a_0,\, a_1,\, a_2,\, a_3)$ denotes the action 
%\[W^{x_0}X^{x_1}Y^{x_2}Z^{x_3} \mapsto (e^{\frac{2\pi i}{n}a_0}W)^{x_0}(e^{\frac{2\pi i}{n}a_1}X)^{x_1}(e^{\frac{2\pi i}{n}a_2}Y)^{x_2}(e^{\frac{2\pi i}{n}a_3}Z)^{x_3}. \]
Let $M$ be a lattice of rank $3$ defined by 
\[
M:=\left\{ (i,j,k,l)\in\mathbb{Z}^4 \, | \, 2i+ 3j+ 5k+ 7l = 0\right\}. 
\]

%It is clear that all monomials of weighted degree $17$ in $\mathbb{P}(2,3,5,7)$ are invariant by the action of $G$. 
%Thus, in this case, $\Delta_{(2,3,5,7), G}$ coincides with $\Delta_{(2,3,5,7)}$. 

By taking a basis $\{ e_1,\, e_2,\, e_3\}$ of $M$ by 
\[
e_1 = (6, 0, -1, -1),\, 
e_2 = (5, -1, 0, -1),\, 
e_3 = (4, -1, -1, 0), 
\]
the polytope $\Delta_{(2,3,5,7)}$ is given by a convex hull of vertices 
\[
\begin{array}{l}
(1, 0, 0),\, (0, 1, 0),\, (0, 0, 1),\, (0, -1, 1),\\
(-1, 1, 0),\, (-1, 2, -1),\, (2, -1, -2),\, (3, -2, -2), 
\end{array}
\]
respectively corresponding to monomials 
\[
\begin{array}{l}
W^7X,\, W^6Y,\, W^5Z,\, XZ^2,\\
 Y^2Z,\, WY^3,\,  X^4Y,\, WX^5. 
\end{array}
\]

The Newton polytope $\Delta_F$ is a convex hull of vertices 
\[
(0, 1, 0),\, (0, 0, 1),\, (0, -1, 1),\, (-1, 1, 0),\, (2, -1, -2).  
\]
Since the polar dual of the face 
\[
{\rm Conv} \left\{ (0, 0, 1),\, (0, -1, 1),\, (2, -1, -2) \right\}
\]
is a rational vertex $(-3/2, 0, -1)$, the Newton polytope is not reflexive. 

Any reflexive polytope $\Delta$ satisfying $\Delta_F\subset\Delta\subset\Delta_{(2,3,5,7)}$ is thus a convex hull of vertices
\[
(1, 0, 0),\, (0, 1, 0),\, (0, 0, 1),\, (0, -1, 1),\, (-1, 1, 0),\, (2, -1, -2),\, (n+1, -n, -n)
\]
with $n=1, 2$, or with $(-1, 2, -1)$. 

As long as there is an edge 
\[
\Gamma = {\rm Conv}\{ (0, -1, 1),\, (-1, 1, 0)\} 
\]
in $\Delta$, since the polar dual of $\Gamma$ is 
\[
\Gamma^* = {\rm Conv}\{ (6, 5, 4),\, (1, 0, -1)\}, 
\]
and thus there are $4$ lattice points in $\Gamma^*$, the polar dual polytope $\Delta':=\Delta^*$ of $\Delta$ should contain an edge with $4$ lattice points. 
However, by a direct check, there does not exist such an edge in or inside of the polytope $\Delta_{(2,3,7,9)}$. 
Thus, the polytope $\Delta'$ is not a subpolytope of $\Delta_{(2,3,7,9)}$. 
Therefore the assertion is verified. 
\QED

\section{Conclusion}
Combining our result with~\cite{MU}, all but singularities that are not assigning an invertible projectivisation, the strange duality for bimodal singularities extends to the polytope duality.

Makiko Mase\\
Universit\"at Mannheim, Lehrstuhl f\"ur Mathematik VI\\
B6, 26, 68131 Mannheim, Germany \\
email: mmase@mail.uni-mannheim.de

\end{document}